\documentclass[a4paper,12pt]{amsart}
\usepackage{fancyheadings,amsmath,amscd,amsthm,amsfonts,latexsym,mathrsfs,amssymb}
\DeclareMathOperator{\grad}{grad}

\DeclareMathOperator{\RicM}{^{\it M}\!Ricci}
\DeclareMathOperator{\RicN}{^{\it N}\!Ricci}

\DeclareMathOperator{\dif}{d}

\renewcommand{\H}{\mathcal{H}}
\newcommand{\V}{\mathcal{V}}
\newcommand{\F}{\mathcal{F}} 
\newcommand{\U}{\mathcal{U}}

\DeclareMathOperator{\dH}{d^{\H}}

\def \G{\Gamma}
\def \l{\lambda}
\def \O{\Omega}
\def \phi{\varphi} 
\def \r{\rho} 

\def \D{\Delta}
\def \Re{\mathbb{R}}
\def \Hq{\mathbb{H}\,}
\def \Co{\mathbb{C}\,}

\begin{document}

\title{A new construction of Einstein self-dual metrics}
\author{Radu Pantilie and John C.\ Wood} 
\thanks{The authors gratefully acknowledge that this work was 
done under E.P.S.R.C. grant number GR/N27897.} 
\email{r.pantilie@leeds.ac.uk\,, j.c.wood@leeds.ac.uk\,.}

\newtheorem{thm}{Theorem}[section]
\newtheorem{lem}[thm]{Lemma}
\newtheorem{cor}[thm]{Corollary}
\newtheorem{prop}[thm]{Proposition}

\theoremstyle{definition}
\newtheorem{defn}[thm]{Definition}
\newtheorem{rem}[thm]{Remark}
\newtheorem{exm}[thm]{Example}

\numberwithin{equation}{section}

\maketitle
\thispagestyle{empty} 
\vspace{-0.5cm}
\section*{Abstract}
\begin{quote}
{\footnotesize
We give a new construction of Ricci-flat self-dual metrics which 
is a natural extension of the Gibbons--Hawking ansatz. We also give characterisations of both these constructions, and explain how they come from harmonic morphisms.} 
\end{quote}

\section*{Introduction}

\indent
In \cite{GibHaw,Haw}\,, G.W.~Gibbons and S.W.~Hawking introduced a construction 
of Einstein self-dual metrics with 
zero scalar curvature (see \cite{LeB} for a thorough discussion 
of this ansatz). The construction is in the spirit of Kaluza--Klein theory with the projection 
of the (local) bundle being a Riemannian submersion followed by a conformal 
transformation.\\ 
\indent
A \emph{harmonic morphism} is a map between Riemannian manifolds which preserves Laplace's 
equation (see Section 2 below). 
In \cite{Bry}\,, R.L.~Bryant gave a local normal form for the metric on 
the domain of a submersive harmonic morphism with one-dimensional fibres. When the 
domain is four-dimensional, 
this local normal form includes that of the Gibbons--Hawking ansatz.\\ 
\indent
In \cite{Pan-thesis,Pan-4to3}\,, it is shown that, from an Einstein four-manifold, there 
are precisely 
three types of harmonic morphism with one-dimensional fibres. The first two types are due 
to R.L.~Bryant and to P.~Baird and J.~Eells, respectively, and lead to the 
Gibbons--Hawking construction and the well-known warped product construction of 
Einstein metrics (see Theorem \ref{thm:GibHaw} and Theorem \ref{thm:warpedprod} below). The third 
type can be seen as a construction of Ricci-flat self-dual metrics.  
We present this new construction in  Section 1 (Theorem \ref{thm:construction}) 
together with the above mentioned 
result from \cite{Pan-thesis,Pan-4to3} reformulated as a classification result for Einstein  four-manifolds whose metric can be written in Bryant's local normal form 
(Theorem \ref{thm:charact}). 
In Section 2 we review some facts on harmonic morphisms. In Section 3 we give proofs of  
the results of Section 1 and show that our construction is a natural extension of the 
Gibbons--Hawking ansatz; indeed all three constructions are characterised by 
an equation \eqref{e:*} which generalizes the monopole equation. We also classify 
the harmonic morphisms with one-dimensional fibres on 
compact Einstein four-manifolds. Our new construction involves solving equation 
\eqref{e:type3} below, which is a particular case of the Beltrami fields equation of 
hydrodynamics (see \cite{KenPlu}\,).  In Section 4 we describe all solutions of \eqref{e:type3}\,, 
both locally and globally, on $S^3$ by a method similar to the one used in \cite{KenPlu} to 
describe solutions of the Beltrami fields equation on $\Re^3$ (see Remark \ref{rem:Beltrami}\,), 
giving all Einstein metrics of the form \eqref{e:construction} below; in fact we show that 
for such metrics the Einstein condition is \emph{equivalent} to the Beltrami fields equation.\\  

\indent 
We are very grateful to G.W.~Gibbons and N.J.~Hitchin for pointing out to us that 
\eqref{e:type3} is a particular case of the Beltrami fields equation and for recommending  
references on this subject. We thank D.M.J.~Calderbank for drawing our attention to his paper 
\cite{Cal}\,, see below, and F.~Belgun for useful comments.

\section{The construction}

\indent
Firstly, we recall the following ansatz of G.W.~Gibbons and S.W.~Hawking 
\cite{GibHaw, Haw} (see also, \cite{LeB}). 

\begin{thm} \label{thm:GibHaw} 
Let $U\subseteq(\Re^3,h)$ be an open subset of the Euclidean three-space. Let $u$ be 
a positive smooth function and $A$ a one-form on $U$. Then, the Riemannian metric on $\Re\times U$ given by  
\begin{equation} \label{e:GibHaw} 
g=u\,h+u^{-1}(\dif\!t+A)^2\qquad\qquad(t\in\Re) 
\end{equation} 
is Einstein if and only if\/ $u$ and $A$ are related by the monopole equation: 
\begin{equation} \label{e:monopole} 
\dif\!u=*\dif\!A   
\end{equation} 
with respect to a suitable choice of orientation on $U$. 
Moreover, if \eqref{e:monopole} holds, then $g$ is Ricci-flat and self-dual. 
\end{thm} 

\indent
Note that \eqref{e:monopole} implies that $u$ is a harmonic function. Moreover, if, 
for example, $U$ is an open ball then, given a 
harmonic function $u$\,,  we can find a one-form $A$ such that \eqref{e:monopole} is satisfied. 
Also, note that $A$ defines a principal connection on $(\Re\times U,U,\Re)$ with curvature form 
$F=\dif\!A$ which is harmonic, i.e., $\dif\!F=0$\,, $\dif^*\!F=0$\,. Conversely, given a harmonic 
two-form $F$ on, say, an open ball $U$ of $\Re^3$\,, solving $\dif\!A=F$ and $\dif\!u=*F$ gives 
a solution to \eqref{e:monopole}\,. 

\indent
Another way of obtaining Einstein metrics, to which our construction is related, is given by the following well-known warped product construction (see \cite[9.109]{Bes}\,). 

\begin{thm} \label{thm:warpedprod} 
Let $(N^n,h)$ be an Einstein manifold of dimension $n$ with $\RicN=c^Nh$\,. Let $g$ be the warped product metric 
on $M^{n+1}=\Re\times N^n$ given by $g=\dif\!t^2+\l^{-2}h$ where $\l:\Re\to(0,\infty)$ is a smooth function.\\ 
\indent
Then, $(M^{n+1},g)$ is Einstein, with $\RicM=c^Mg$\,, if and only if $\l$ satisfies the following equation:  
\begin{equation} \label{e:oden} 
\frac{c^M}{n}\,\l^2-\frac{c^N}{n-1}\,\l^4+\bigl(\l'\bigr)^2=0\;.  
\end{equation} 
Moreover, if $(N^n,h)$ has constant curvature  
(note that this is automatically satisfied if $\dim N=3$), then $(M^{n+1},g)$ also 
has constant curvature. 
\end{thm}  

\indent
Next, we give the new construction based on the \emph{Beltrami fields equation},
see Proposition \ref{prop:solv} for solutions. 

\begin{thm} \label{thm:construction} 
Let $U\subseteq(S^3,h)$ be an open subset of the three-sphere endowed 
with its canonical metric.\\ 
\indent
Let $A$ be a one-form on $U$. Then, the Riemannian metric on 
$(0,\infty)\times U$ given by 
\begin{equation} \label{e:construction} 
g=\r^2h+\r^{-2}(\r\dif\!\r+A)^2\qquad\qquad(\r\in(0,\infty)\,) 
\end{equation} 
is Einstein if and only if the following \emph{Beltrami fields equation} holds on $U$:  
\begin{equation} \label{e:type3} 
\dif\!A+2*\!A=0
\end{equation} 
with respect to a suitable choice of orientation on $U$.  
Moreover, if \eqref{e:type3} holds, then $g$ is Ricci-flat and self-dual.  
\end{thm} 

\begin{rem} 
D.M.J.~Calderbank \cite[Theorem V]{Cal} reduces the problem of finding scalar flat
K\"ahler self-dual $4$-manifolds 
to certain affine monopole equations on an Einstein-Weyl 3-manifold. For suitable 
choice of gauge, these equations reduce to the Beltrami fields equation.   
\end{rem}

\indent
The proof of Theorem \ref{thm:construction} will be given in Section 3\,.\\ 

\indent
By a \emph{local principal bundle} $(M^{n+1},N^n,S^1)$ (with \emph{principal connection}  
$\H\subseteq TM$\,) we shall mean the restriction 
of a principal bundle $(\widetilde{M}^{n+1},\widetilde{N}^n,S^1)$ 
(with principal connection $\widetilde{\H}\subseteq T\widetilde{M}$\,) to an open subset 
$M^{n+1}$ of $\widetilde{M}^{n+1}$\,. 
Then, the above three constructions admit the following characterisation which will 
be proved in Section 3. 

\begin{thm} \label{thm:charact} 
Let $(N^3,h)$ be a Riemannian three-manifold and let $(M^4,N^3,S^1)$ be a
local principal bundle  
endowed with a principal connection $\H\subseteq TM$\,. Define a Riemannian metric 
$g$ on $M^4$ by 
\begin{equation} \label{e:Bryant} 
g=\l^{-2}\,\phi^*(h)+\l^2\,\theta^2 
\end{equation} 
where $\l$ is a positive smooth function on $M^4$\,, $\phi:M^4\to N^3$ is the projection 
of the local principal bundle $(M^4,N^3,S^1)$ and $\theta$ is the connection form of $\H$\,.\\ 
\indent
Suppose that $(M^4,g)$ is Einstein.\\ 
\indent
{\rm (a)} Then one of the following assertions holds:\\ 
\indent
\quad{\rm (i)} $S^1$ acts by isometries on $M^4$ (in particular, the fibres of $\phi$ form 
a Riemannian foliation);\\ 
\indent
\quad{\rm (ii)} $(M^4,g)$ and $(N^3,h)$ have constant sectional curvature and, locally, $g$ is given by Theorem \ref{thm:warpedprod}\,, with $n=3$\,;\\ 
\indent
\quad{\rm (iii)} $(N^3,h)$ has positive constant sectional curvature and, up to homotheties, 
$g$ is locally given by Theorem \ref{thm:construction}\,.\\ 
\indent
{\rm (b)} Furthermore, if $(N^3,h)$ has constant curvature then {\rm (i)} can be replaced by\\ 
\indent
\quad{\rm (i$'$)} $(N^3,h)$ is flat and $g$ is locally given by Theorem \ref{thm:GibHaw}\,.\\  
\indent
Thus if $(N^3,h)$ has constant curvature, the ansatz \eqref{e:Bryant} gives an Einstein metric if and only if it is one of the three constructions: Gibbons--Hawking (Theorem \ref{thm:GibHaw}),  
warped product (Theorem \ref{thm:warpedprod}) or Beltrami fields (Theorem \ref{thm:construction}). 
\end{thm}

\begin{rem} 
1) If both (i) and (ii) of Theorem \ref{thm:charact} occur then, locally, $(M^4,g)$ is the 
Riemannian product $\Re\times(N^3,h)$\,; whilst, if both (ii) and (iii) occur, then, 
locally and up to homotheties, $(N^3,h)$ is the three-sphere, $(M^4,g)$ is the Euclidean four-space 
and \eqref{e:Bryant} corresponds to the polar coordinates. However, (i) and (iii) cannot 
occur simultaneously.\\  
\indent
2) Note that a map $\phi:(M^{n+1},g)\to(N^n,h)$ is a submersive harmonic morphism 
if and only if it is locally as in the hypothesis of Theorem \ref{thm:charact} 
(see Theorem \ref{thm:Bry1} below).\\ 
\indent 
3) A similar result (\cite{PanWood-d}) can be given for a Riemannian manifold $N$ of any 
dimension $n\geq4$\,. In 
that case, $g=\l^{-2}\,\phi^*(h)+\l^{2n-4}\,\theta^2$ (cf.\ \cite{Bry}), 
and, (i$'$) and (iii) can only occur in the trivial cases $u={\rm constant}$ and $A=0$\,.\\ 
\indent 
4) Assertion (i) of Theorem \ref{thm:charact} is equivalent to the assertion that 
$\l$ be constant along the fibres of\/ $\phi$\,, whereas, in case (ii)\,, $\l$ is constant along 
horizontal curves (i.e., tangent to $\H$), see, also, Proposition \ref{prop:BaiEel} and 
Proposition \ref{prop:Bry}, below.\\ 
\indent
5) In assertion (iii) of Theorem \ref{thm:charact}\,, we have 
$\rho=\l^{-1}$ where $\rho$ is as in Theorem \ref{thm:construction}\,. 
Similarly, in assertion (i$'$) of Theorem \ref{thm:charact}\,, we have $u=\l^{-2}$ 
where $u$ is as in Theorem \ref{thm:GibHaw}\,. 
\end{rem}

\section{Some facts on harmonic morphisms}

\indent
Both Theorem \ref{thm:construction} and Theorem \ref{thm:charact} have their origin  
in the theory of harmonic morphisms; in this section, we recall some basic 
facts from that theory (see \cite{BaiWoo2} for a general account and \cite{Gudbib} for a 
frequently updated bibliography). For simplicity, 
\emph{from now on all the manifolds are assumed to be smooth, connected and orientable}. 

\begin{defn} 
A \emph{harmonic morphism} is a smooth map $\phi:(M^m,g)\to(N^n,h)$ between Riemannian manifolds 
which pulls back harmonic functions to harmonic functions, i.e., for any harmonic 
function $f:(U,h|_h)\to\Re$ defined on an open subset $U$ of $N^n$\,, with 
$\phi^{-1}(U)$ non-empty, $f\circ\phi:(\phi^{-1}(U),g|_{\phi^{-1}(U)})\to\Re$ is a 
harmonic function. 
\end{defn} 

\indent
To state the basic characterisation result for harmonic morphisms we also need the 
following.  

\begin{defn} 
A smooth map $\phi:(M^m,g)\to(N^n,h)$ between Riemannian manifolds is \emph{horizontally 
(weakly) conformal} if, at each point $x\in M$\,, \emph{either} $\dif\!\phi_x=0$\,, \emph{or}  
$\dif\!\phi_x:T_xM\to T_{\phi(x)}N$ is surjective and its restriction to the 
horizontal space $\H_x=({\rm ker}\dif\!\phi_x)^{\perp}$ is a conformal (linear) map 
$(\H_x,g_x|_{\H_x})\to(T_{\phi(x)}N,h_{\phi(x)})$\,. Denote the conformality factor by $\l(x)$\,. 
The resulting function $\l$ is called the \emph{dilation} of $\phi$\,. The dilation 
is smooth outside the set of critical points and can be extended to a continuous function 
on $M^m$\,, with $\l^2$ smooth, by setting it equal to zero on the set of 
critical points.\\ 
\indent
A smooth map is called \emph{horizontally homothetic} if it is horizontally conformal 
with dilation constant along horizontal curves. 
\end{defn} 

\indent
The following theorem is due to B.~Fuglede \cite{Fug} and T.~Ishihara \cite{Ish}\,. 

\begin{thm} \label{thm:FugIsh} 
A smooth map between Riemannian manifolds is a harmonic morphism if and only if 
it is a harmonic map which is horizontally weakly conformal. 
\end{thm}

\indent
\emph{Any Riemannian submersion with minimal fibres is a harmonic morphism.} 
This can be seen either directly or by applying Theorem \ref{thm:FugIsh}\,. 
Other simple examples 
are \emph{radial projection} $\Re^{n+1}\setminus\{0\}\to S^n$ ($n\geq1$) defined by 
$x\mapsto x/|x|$ and the \emph{Hopf polynomial} $\Co^2\to\Re^3$ defined by 
$(z_1,z_2)\mapsto\bigl(|z_1|^2-|z_2|^2,2z_1\overline{z_2}\bigr)$\,. 
Note that radial projection 
is a horizontally homothetic submersion with geodesic fibres, whilst the fibres 
of the Hopf polynomial are tangent to a Killing vector field. Both of these examples 
are particular cases of more general constructions which we now present. 

\begin{prop}[\cite{BaiEel}] \label{prop:BaiEel} 
Let $\phi:(M^m,g)\to(N^n,h)$ be a non-constant horizontally weakly conformal map. If $\dim N=2$\,, 
then $\phi$ is a harmonic morphism if and only if its fibres are minimal at regular points. 
If\/ $\dim N\neq2$\,, then any two of the 
following assertions imply the third:\\ 
\indent
{\rm (i)} $\phi$ is a harmonic morphism;\\  
\indent
{\rm (ii)} the fibres of $\phi$ are minimal at regular points;\\ 
\indent
{\rm (iii)} $\phi$ is horizontally homothetic. 
\end{prop}  

\indent
{}From Proposition \ref{prop:BaiEel} it follows that \emph{any horizontally homothetic submersion 
with minimal fibres is a harmonic morphism}. 

\begin{prop}[\cite{Bry}] \label{prop:Bry} 
For $n\geq3$\,, let $\phi:(M^{n+1},g)\to N^n$ be a surjective submersion with 
one-dimensional fibres which form a Riemannian foliation. Then the following  
assertions are equivalent:\\ 
\indent
{\rm (i)} there exists a Riemannian metric $h$ on $N$ with respect to which 
the map $\phi:(M^{n+1},g)\to(N^n,h)$ is a harmonic morphism;\\  
\indent
{\rm (ii)} there exists a non-zero Killing vector field tangent to the fibres 
of $\phi$\,. 
\end{prop} 

\indent
{}From \cite{Bry}\,, we also recall the following characterisation. 

\begin{thm} \label{thm:Bry1} 
Let $(M^{n+1},N^n,S^1)$ be a principal bundle with projection 
$\phi:M^{n+1}\to N^n$ and endowed 
with a principal connection $\H\subseteq TM$\,. Let $h$ be a Riemannian metric on $N^n$ and 
$\l$ a smooth positive function on $M^{n+1}$\,.\\ 
\indent
Define a Riemannian metric on $M^{n+1}$ by 
\begin{equation} \label{e:Bry} 
g=\l^{-2}\,\phi^*(h)+\l^{2n-4}\,\theta^2 
\end{equation} 
where $\theta$ is the connection form of $\H$\,. Then, 
$\phi:(M^{n+1},g)\to(N^n,h)$ is a harmonic morphism.\\  
\indent
Conversely, any submersive harmonic morphism with one-dimensional fibres is \emph{locally} 
of this form, up to isometries. 
\end{thm} 

\indent
See \cite{Pan}\,,\,\cite{Pan-thesis} for a proof of Theorem \ref{thm:Bry1} and a more 
explicit version of the converse.   

\begin{rem} \label{rem:harmorphs} 
1) In the notations of Theorem \ref{thm:Bry1}\,, let $V$ be the vertical vector field with  $\theta(V)=1$\,. Obviously, $V$ is the infinitesimal generator of the $S^1$ action and $g(V,V)=\l^{2n-4}$\,. We call $V$  \emph{the fundamental (vertical) vector field}.\\ 
\indent
2) If $(M^4,g)$ is as in Theorem \ref{thm:GibHaw} (respectively, Theorem \ref{thm:warpedprod}\,,  Theorem \ref{thm:construction}\,) and $N^3$ is an 
open subset of $\Re^3$ (respectively, a constant curvature three-manifold, $S^3$\,),   
then the canonical projection $\phi:(M^4,g)\to(N^3,h)$ is a harmonic morphism. 
Furthermore, if $g$ is given by the Gibbons--Hawking construction 
then, obviously, its fibres are generated by a Killing vector field, whilst if $g$ 
is given by Theorem \ref{thm:construction}\,, with $A\neq0$\,, then neither is $\phi$  
horizontally homothetic nor do its fibres form a Riemannian foliation.\\ 
\indent
3) Recall \cite{Bai} that any harmonic morphism with one-dimensional fibres from a Riemannian 
manifold of dimension at least five is submersive, whilst if the domain has dimension four 
then the set of critical points is discrete. 
\end{rem}

\section{Characterisations of the construction and some related results}

\indent
{}From the previous section it follows that classifying Einstein metrics which can be 
locally written in the form \eqref{e:Bry} is the same as classifying harmonic morphisms with 
one-dimensional fibres from Einstein manifolds. Therefore the results of this section will be 
given in the language of harmonic morphisms. The reader who is not primarily interested 
in harmonic morphisms can easily rewrite all these results in the language of Section 1.\\

\begin{prop} \label{prop:*} 
Let $(M^4,g)$ be an Einstein four-manifold and $\phi:(M^4,g)\to(N^3,h)$ a submersive 
harmonic morphism to a Riemannian three-manifold.\\ 
\indent
Then, the following assertions are equivalent:\\ 
\indent
{\rm (i)} $(N^3,h)$ has constant curvature;\\  
\indent
{\rm (ii)} the following equation holds:  
\begin{equation} \label{e:*} 
\dH(\l^{-2})=*_\H\O\;;
\end{equation} 
here\/ $\dH$ is the differential followed by the orthogonal projection onto $\H$\,, $\l$ is the dilation 
of\/ $\phi$\,, $*_\H$ is the Hodge star operator on\/ $\bigl(\H,\phi^*(h)\bigr)$ with respect to a 
suitably chosen orientation,  and\/ $\O$ is the curvature form of the horizontal 
distribution (i.e., in the notation of Theorem \ref{thm:Bry1}\,, $\O=\dif\theta$). 
\end{prop} 
\begin{proof} 
The following formula is a consequence of \cite[(B.0.23)]{Pan-thesis}\,.  
\begin{equation} \label{e:riccixy*} 
\begin{split} 
\RicM|_{\H\otimes\H}=\phi^*\bigl(\!&\RicN\bigr)-\l^{-2}\bigl(\,\D\!^M(\log\l)
+\tfrac14\,\l^6\,|\O|_h^2\,\bigr)\,\phi^*(h)\\ &+\tfrac12\,\l^4\,\bigl(*_\H\O\bigr)\otimes\bigl(*_\H\O\bigr)-2\dH(\log\l)\otimes\dH(\log\l)\,. 
\end{split} 
\end{equation} 
\indent
Assume (ii)\,. Then from \eqref{e:riccixy*} and a Schur-type lemma (see \cite{Bes}\,) it 
easily follows that $(N^3,h)$ is Einstein and hence of constant curvature.\\ 
\indent
Conversely, suppose that (i) holds. Recall (see 
\cite[Theorem 5.26]{Bes}\,) that any 
Einstein manifold can be given a real-analytic structure. Then, as in 
\cite[Proposition 1.4]{PanWood-d}\,, it follows that all the objects appearing in  \eqref{e:riccixy*} are real-analytic.\\ 
\indent
{}From \eqref{e:riccixy*} it follows that $\O=0$ if and only if $\dH\!\l=0$\,. 
Otherwise, by real-analyticity, both are non-zero on a dense open subset of $M^4$ and 
we may choose 
a real-analytic positive orthonormal local frame $\{X_1,\,X_2,\,X_3\}$ for 
$(\H,\phi^*(h))$ such that $X_i(\l)$ is nowhere zero for all $i=1,2,3$\,. 
{}From  \eqref{e:riccixy*} it follows that 
\begin{equation} \label{e:before*} 
\bigl(*_\H\O\bigr)(X_i)\,\bigl(*_\H\O\bigr)(X_j)=X_i(\l^{-2})\,X_j(\l^{-2}) 
\end{equation}
for any $i,\,j=1,2,3$\,, $i\neq j$\,. Equation \eqref{e:*} follows easily from \eqref{e:before*}\,. 
\end{proof} 

\indent
Next, we give the proof of Theorem \ref{thm:charact}. 

\begin{proof}[Proof of Theorem \ref{thm:charact}] 
Part (a) follows from the main result of \cite{Pan-4to3} (see 
\cite[Corollary 3.4.5]{Pan-thesis}\,, and the proof of \cite[Theorem 2.8]{PanWood-d}\,),  
by noting that Theorem \ref{thm:charact} can be viewed as a classification result for 
harmonic morphisms with one-dimensional fibres on Einstein four-manifolds.\\ 
\indent
Part (b) is a consequence of Proposition \ref{prop:*}\,. Indeed, suppose that both $(M^4,g)$ 
and $(N^3,h)$ are Einstein. If $\O=0$ on $M$ then, from \eqref{e:*}\,, it follows that $\l$  
is constant along horizontal curves, and Proposition \ref{prop:BaiEel} implies  
that, locally, $g$ is a warped product. Then, \eqref{e:oden} (with $n=3$) follows, 
for example, from  
\cite[9.109]{Bes} (see also \cite{Pan-thesis}), and hence assertion (ii) of 
Theorem \ref{thm:charact} holds.\\ 
\indent
If $\O\neq0$\,, then, as in the previous proof, $\O$ is real-analytic; hence it is non-zero 
on a dense open subset of $M^4$\,. Now, note that the right hand side of \eqref{e:*} is basic.  
Hence $V(X(\l^{-2}))=0$ for any basic vector field $X\in\G(\H)$ where $V$ is the fundamental 
vector field (Remark \ref{rem:harmorphs}(1)\,). 
But $V$ commutes with basic vector fields (because $\H$ is a local principal connection on $\phi$),  
and hence $V(\l^{-2})$ is constant along horizontal curves. It follows that, if $V(\l^{-2})$ is 
non-constant then $\H$ is integrable, equivalently, $\O=0$\,. Thus we must have that 
$V(\l^{-2})=c$ for some constant $c\in\Re$\,.\\ 
\indent
If $c=0$\,, then $\l$ is constant along the 
fibres of $\phi$ (equivalently, $g$ is as in assertion (i) of Theorem \ref{thm:charact}\,), 
and from \eqref{e:*} it will follow that $g$ is locally given by 
the Gibbons--Hawking construction, once we have shown that $(N^3,h)$ is flat; this 
will follow from \eqref{e:ricci*} below.\\ 
\indent
If $c\neq0$\,, then $(1/c)\dif(\l^{-2})$ is a (flat) principal connection on $\phi$\,. 
Let $A\in\G(T^*N)$ be a local connection form of $\H$ with respect to 
$(1/c)\dif(\l^{-2})$\,, that is, $A$ is the one-form on $N$ which satisfies 
\begin{equation} \label{e:locconn} 
\theta=\frac1c\,\dif(\l^{-2})+\phi^*(A)\;. 
\end{equation} 
\indent
{}From \eqref{e:*} and \eqref{e:locconn}\,, it follows that 
\begin{equation} 
-c\,\phi^*(A)=\dH(\l^{-2})=*_\H\,\O=*_\H\;\phi^*(\dif\!A)=\phi^*(*\dif\!A)\;.   
\end{equation} 
Hence $\dif\!A+c*A=0$ which implies assertion (iii)\,, except for the fact that $N^3$ has constant 
sectional curvature equal to $c^2/4$ which we shall now prove.\\ 
\indent
Let $\phi:(M^4,g)\to(N^3,h)$ be a submersive harmonic morphism between Riemannian manifolds 
of dimension four and three, respectively. If \eqref{e:*} is satisfied then, by applying 
\cite[Lemma B.0.19]{Pan-thesis}\,, we get that the Ricci tensors of $(M,g)$ and $(N,h)$ 
satisfy the following relations: 
\begin{equation} \label{e:ricci*} 
\begin{split}  
&\RicM|_{\V\otimes\V}=0\;,\quad\RicM|_{\V\otimes\H}=0\;,\\ 
&\RicM|_{\H\otimes\H}=\phi^*\bigl(\RicN\bigr)-\frac{c^2}{2}\,\phi^*(h)\;. 
\end{split} 
\end{equation} 
\indent
{}From \eqref{e:ricci*} it follows easily that, if $(M^4,g)$ is Einstein and \eqref{e:*} 
holds, then it is Ricci-flat and $(N^3,h)$ has constant sectional curvature equal to $c^2/4$\,. 
\end{proof} 

\begin{rem} 
The proof of Theorem \ref{thm:charact} shows that \eqref{e:*} unifies the $S^1$-monopole 
equation and Beltrami fields equation. Therefore our construction 
may be considered to be a natural extension of the Gibbons--Hawking ansatz. 
\end{rem} 

\indent
Theorem \ref{thm:construction} follows from the following result. 

\begin{thm} \label{thm:constr} 
Let $\phi:(M^4,g)\to(N^3,h)$ be a surjective submersive harmonic morphism such that 
$\dH(\l^{-2})=*_\H\O$ with $V(\l^{-2})=c\,(\in\Re)$\,.\\ 
\indent
Then the following assertions are equivalent:\\ 
\indent
{\rm (i)} $(M^4,g)$ is Einstein;\\ 
\indent
{\rm (ii)} $(N^3,h)$ has constant sectional curvature equal to $c^2/4$\,.\\ 
\indent
Moreover, if\/ {\rm (i)} or\/ {\rm (ii)} holds, then $(M^4,g)$ is Ricci-flat and 
self-dual. 
\end{thm} 
\begin{proof} 
The equivalence of (i) and (ii) follows from \eqref{e:ricci*}\,.\\ 
\indent
To establish the last statement, we may assume for simplicity that, if $c\neq0$\,, $(N^3,h)$ 
is the three-sphere 
of radius $2/c$ with its canonical metric and, if $c=0$\,, $(N^3,h)$ is the Euclidean three-space. 
Take $p$ to be the Hopf fibration $S^3\to S^2$ if $c\neq0$\,, or an orthogonal 
projection $\Re^3\to\Re^2$ if $c=0$\,.\\
\indent
Then, $\psi=p\circ\phi$ is a submersive harmonic morphism with two-dimensional fibres. 
Let $\F={\rm ker}\dif\!\psi$ and define $J$ to be the (negatively oriented) almost 
Hermitian structure on $(M^4,g)$ given by rotation through angle $-\pi/2$ on $\F$ and 
rotation through angle $\pi/2$ on  $\F^{\perp}$\,.\\ 
\indent
A straightforward calculation shows that $(M^4,J,g)$ is a K\"ahler manifold if and only 
if $\dH(\l^{-2})=*_\H\O$ with $V(\l^{-2})=c$\,. It follows that $(M^4,g)$ is 
hyper-K\"ahler and, in particular, Ricci-flat and self-dual (see \cite{Bes}\,). 
(Note that the orientation of $M^4$ is given as follows: orient $\V={\rm ker}\dif\!\phi$ 
by $\theta$ and then choose the orientation on $M^4$ such that the canonical vector bundle isomorphisms $TM=\V\oplus\H$ and $\H=\phi^*(TN)$ are orientation preserving.) 
\end{proof} 

\begin{rem} 
1) The second part of the proof of Theorem \ref{thm:constr} was inspired by C.~LeBrun's  discussion of the Gibbons--Hawking ansatz \cite{LeB} (see also \cite{Woo-4d}\,, \cite[Theorem 3.5]{CheDon}).\\ 
\indent
2) In \cite{PanWoo-sd}\,, the method is developed to give the following:\\ 
\indent
\emph{Let\/ $(N^3,h)$ be a constant curvature three-manifold and let\/ $A$ be a 
real-analytic one-form on $N$. Define a Riemannian metric on\/ $(0,\infty)\times N^3$ by 
\begin{equation} \label{e:1sdconstruction} 
g=\rho\,h+\rho^{-1}(\dif\!\rho+A)^2\qquad\qquad(\rho\in(0,\infty)\,)\;. 
\end{equation} 
\indent
Then\/ $g$ is self-dual (respectively, anti-self-dual) if and only if the following Beltrami 
fields equation holds on\/ $N^3$:}
\begin{equation} \label{e:half} 
\dif\!A=-*A\qquad\textit{(respectively,\:\,$\dif\!A=*A$\,)}\;.
\end{equation} 
\indent
Note that $(N^3,h)$ may have constant sectional curvature of sign unrelated to that in 
\eqref{e:half} and that $g$ is not, in general, Einstein.\\ 
\indent
We shall give details on this construction in \cite{PanWoo-sd}\,, 
together with another new construction of self-dual metrics. 
\end{rem} 

\indent
We end this section with the classification of harmonic morphisms with one-dimensional fibres 
on compact Einstein four-manifolds. For this we need the following, which improves one of the 
statements of \cite[Theorem 2.9]{PanWoo-Vran}\,.

\begin{prop} \label{prop:critp-Ricci-flat}  
Let $(M^4,g)$ be an Einstein four-manifold and let $\phi:(M^4,g)\to(N^3,h)$ be a harmonic 
morphism with one-dimensional fibres to a Riemannian three-manifold.\\ 
\indent
If $\phi$ has critical points, then $(M^4,g)$ is Ricci-flat and there exists a (real-analytic)   Killing vector field tangent to the fibres of $\phi$ whose zero set is equal to the set 
of critical points of $\phi$\,.  
\end{prop} 
\begin{proof} 
By a result of P.~Baird \cite{Bai}\,, a harmonic morphism $\phi:(M,g)\to(N,h)$ with 
one-dimensional fibres is submersive if $\dim M>4$\,, whilst, if $\dim M=4$\,, the 
set of critical points of $\phi$ is discrete.\\ 
\indent
Now Theorem \ref{thm:charact} is equivalent to a classification result for harmonic 
morphisms with one-dimensional fibres on Einstein four-manifolds 
\cite[Corollary 1.9]{Pan-4to3}\,,\,\cite[Corollary 3.4.5]{Pan-thesis}\,. As already  explained, 
if assertion (ii) of Theorem \ref{thm:charact} holds, then $\phi$ is a harmonic map which 
is horizontally homothetic, and hence, by a result of B.~Fuglede \cite{Fugcrit}\,, it is 
submersive.\\ 
\indent
If assertion (iii) of Theorem \ref{thm:charact} holds then $V(\l^{-2})$ is a \emph{non-zero} 
constant where, as before, $\l$ is the dilation of $\phi$ and $V$ is the fundamental vector field 
(Remark \ref{rem:harmorphs}(1)\,).  Because the critical points of $\phi$ are isolated, 
from \cite{ChuTim} it follows that, in the neighbourhood of a critical point, $\phi$ is 
topologically equivalent to the cone over the Hopf fibration $S^3\to S^2$. Therefore there 
exists a connected component of a fibre of $\phi$ which is diffeomorphic to $S^1$, and hence 
at some point we must have $V(\l^{-2})=0$\,. Thus, if (iii) of Theorem \ref{thm:charact} 
holds, then $\phi$ is submersive.\\ 
\indent
We have thus shown that, if $\phi$ has critical points, then its fibres are generated by a 
Killing vector field, namely $V$\,, which can be extended to a real-analytic Killing vector 
field on $(M^4,g)$ by setting it equal to zero at the critical points 
(see the proof of \cite[Corollary 3.3]{Pan-4to3} or \cite[Corollary 3.6.3]{Pan-thesis}\,).\\   
\indent
Because $\l$ and $\O$ are basic they are locally the pull-backs of a function and a two-form,  respectively, 
which are defined on $N^3$\,. For simplicity, we shall denote the corresponding objects 
on $N^3$ by the same letters $\l$ and $\O$. Recall that $\RicM=c^Mg$ and  
$g|_\H=\l^{-2}\phi^*(h)$\,. Thus \eqref{e:riccixy*} can be written as an equation on $N^3$\,. 
Furthermore, by applying \cite[(B.0.25)]{Pan-thesis}\,, the corresponding equation on $N^3$ 
can be written as follows: 
\begin{equation} \label{e:riccixy*N} 
\RicN=2c^M\l^{-2}\,h-\tfrac12\l^4(*\O)\otimes(*\O)+2\l^{-2}\dif\!\l\otimes\dif\!\l\;. 
\end{equation} 
\indent
{}From \eqref{e:riccixy*N}, it follows that any vector orthogonal to both 
$(*\O)^{\sharp}$ and $\grad\l$ 
is an eigenvector for $\RicN$\,, the corresponding eigenvalue being $c^M\l^{-2}$\,. 
Hence, if $c^M\neq0$\,, this eigenvalue tends to $\infty$ as we approach a critical value 
of $\phi$\,, which is obviously impossible (apply, for example, \cite[Lemma 2.1]{PanWood-d}\,)\,. 
Thus $c^M=0$\,, i.e., $\RicM=0$\,. 
\end{proof}  

\indent
By applying Proposition \ref{prop:critp-Ricci-flat} we obtain the following result, which improves 
\cite[Theorem 3.8\,,\,Theorem 4.11]{Pan-4to3} (see 
\cite[Theorem 4.1(iii)\,,\,Theorem 4.5]{PanWoo-toprestr}\,). 

\begin{thm} \label{thm:compact4d} 
Let $(M^4,g)$ be a compact Einstein four-manifold and $\phi:(M^4,g)\to(N^3,h)$ a non-constant 
harmonic morphism to a Riemannian three-manifold.\\ 
\indent
Then, up to homotheties and Riemannian coverings, $\phi$ is the canonical projection 
$T^4\to T^3$ between flat tori. 
\end{thm} 
\begin{proof} 
If $\phi$ is submersive this follows from \cite[Theorem 4.1]{PanWoo-toprestr}\,.\\ 
\indent
If $\phi$ has critical points then, by Proposition \ref{prop:critp-Ricci-flat}\,,  
$(M^4,g)$ is Ricci-flat and there exists a Killing vector 
field $V$ tangent to the fibres of $\phi$ whose zero set is equal to the set of critical points 
of $\phi$\,. But, since $(M^4,g)$ is compact and Ricci-flat, by a well-known result of 
S.~Bochner (see \cite[1.84]{Bes}\,) $V$ is actually parallel and hence of constant norm. 
It follows that $\phi$ must be submersive and the theorem is proved. 
\end{proof}

\section{Solving $\dif\!A+2*\!A=0$ on $S^3$}

\indent 
In what follows we regard $S^3$ as the Lie group Sp($1$) of unit quaternions. Let $\{X_j\}$ 
be an orthonormal positively oriented left invariant frame on $S^3$\,, 
$\{\theta_j\}$ its dual and $\pi:S^3\to S^2(\frac12)$ the Hopf fibration chosen so that its 
fibres are tangent to $X_3$\,; here $S^2(\tfrac12)$ denotes the sphere of radius $1/2$\,. 
For any three-dimensional submanifold (with boundary) $N\subseteq S^3$ 
we shall denote the space of harmonic function on $N$ by $\H(N)$\,,    
the space of harmonic functions on $N$ which are basic 
with respect to $\pi$ by $\H_b(N)$\,, and the space of harmonic one-forms on $N$ by $\H^1(N)$\,.

\begin{prop} \label{prop:solv} 
There exists a subbase $\U$ for the topology of $S^3$ such that, for each $U\in\U$\,, the closure 
$\overline{U}$ is a three-dimensional submanifold with boundary, and we have an isomorphism 
\begin{equation} \label{e:isomorph} 
\begin{split} 
\bigl\{&A\in\G(T^*\overline{U})\,|\,\dif\!A+2*\!A=0\:{\rm on}\:\overline{U}\,\bigr\} 
\overset{\sim}{\longrightarrow}\frac{\H(\overline{U})\times\H(\overline{U})\times\H_b(\overline{U})} 
{\H^1(\overline{U})}\\  
\end{split} \end{equation} 
given by 
\begin{equation*} 
\begin{split} 
&A=A_j\,\theta^j\longmapsto\bigl(A_1\,,A_2\,,A_3+\psi_{\overline{U}}(X_1(A_1)+X_2(A_2))\bigr) 
\quad\bigl({\rm mod}\,\H^1(\overline{U})\bigr)  
\end{split} 
\end{equation*} 
where $\psi_{\overline{U}}:\H(\overline{U})\to\H(\overline{U})$ is a certain injective 
linear map such that $X_3\bigl(\psi_{\overline{U}}(u)\bigr)=u$ for any 
$u\in\H(\overline{U})$\,. 
\end{prop} 
\begin{proof} 
Firstly, we construct the subbase $\U$ and the functor $\psi$\,. For each closed $2$-ball $D$ in 
$S^2$ let $s:D\to S^3$ be a section of $\pi$\,. Because $\pi^{-1}(D)$ is a solid torus we 
can find an open $3$-ball $U$ whose closure $\overline{U}\subseteq\pi^{-1}(D)$ is 
a closed $3$-ball such that 
$s(D)\subseteq\overline{U}$ 
and $\pi|_{\overline{U}}$ has connected fibres. We define $\U$ to be formed of all such $U$.\\ 
\indent
Let $U\in\U$\,. To define $\psi_{\overline{U}}$ let $u\in\H(\overline{U})$ and take 
$w:\overline{U}\to\Re$ to be the smooth function characterised by 
\begin{equation} \label{e:w} 
X_3(w)=u\;,\quad w|_{s(\pi(\overline{U}))}=0\;. 
\end{equation} 
Because $X_3$ is a Killing vector field, we have 
\begin{equation*} 
X_3(\D w)=\D(X_3(w))=\D u=0\:. 
\end{equation*} 
Hence we can find $\widetilde{w}:\pi(\overline{U})\to\Re$ such that 
\begin{equation} \label{e:Dwbasic}
\pi^*(\widetilde{w})=\D w\;.  
\end{equation} 
Take $w_1:\pi(\overline{U})\to\Re$ to be the unique solution 
of the following Dirichlet problem: 
\begin{equation} \label{e:Dirichlet} 
\D\!^{S^2}w_1=\widetilde{w}\;,\quad w_1|_{\partial\pi(\overline{U})}=0\;. 
\end{equation} 
\indent
Define $\psi_{\overline{U}}(u)=w-\pi^*(w_1)$\,; we show that  $X_3\bigl(\psi_{\overline{U}}(u)\bigr)=u$ and $\D(\psi_{\overline{U}}(u))=0$\,. 
Indeed, by \eqref{e:w}\,, we have 
\begin{equation*} 
X_3\bigl(\psi_{\overline{U}}(u)\bigr)=X_3(w)-X_3(\pi^*(w_1))=u-0=u\;. 
\end{equation*} 
Also, by \eqref{e:Dwbasic} and \eqref{e:Dirichlet}\,, and since $\pi$  intertwines Laplacians 
(see, for example, \cite[Chapter 4]{BaiWoo2}\,) we have 
\begin{equation*} 
\begin{split} 
\D(\psi_{\overline{U}}(u))&=\D w-\D(\pi^*(w_1))\\ 
	&=\pi^*(\widetilde{w})-\pi^*(\D\!^{S^2}w_1)\\ 
	&=\pi^*(\widetilde{w}-\D\!^{S^2}w_1)=0\;. 
\end{split} 
\end{equation*} 
\indent
For a one-form $A=A_j\,\theta^j$ on $S^3$ set $\D A=(\D A_j)\,\theta^j$\,. 
(Note that $\D$ is \emph{not} the Hodge Laplacian $\dif^*\!\dif+\dif\dif^*$\,.)
To prove that \eqref{e:isomorph} is an isomorphism, it is sufficient to prove the 
following two facts:\\  
\indent
(i) The map $\H(\overline{U})\times\H(\overline{U})\times\H_b(\overline{U})
\longrightarrow\bigl\{A\,\big|\,\D A=0\,,\:\dif^*\!A=0\,\bigr\}$ given by 
\begin{equation*}  
(u_1\,,u_2\,,u_3)\longmapsto\bigl(u_1\,,u_2\,,u_3-\psi_{\overline{U}}(X_1(u_1)+X_2(u_2))\bigr)\:\;     {\rm is\:\; an\:\; isomorphism;} 
\end{equation*} 
\indent   
(ii) $\bigl\{A\,\big|\,\D A=0\,,\:\dif^*\!A=0\,\bigr\}=\bigl\{A\,\big|\,\bigl((*\dif)^2+2(*\dif)\bigr)(A)=0\,, 
\:\dif^*\!A=0\,\bigr\}$\,.\\ 
\indent
The proof of (i) follows easily from the properties of $\psi_{\overline{U}}$ and the fact that 
$\dif^*\!A=-\bigl(X_1(A_1)+X_2(A_2)+X_3(A_3)\bigr)$\,.\\ 
\indent
To prove (ii) we use the Dirac operator ${\rm D}:\G(S^3\times\Hq)\to\G(S^3\times\Hq)$ 
defined by 
\begin{gather*} 
{\rm D}= \begin{pmatrix} 
\, 0  & -X_1 & -X_2 & -X_3 \\ 
\,X_1 &  0   & -X_3 &  X_2 \\ 
\,X_2 &  X_3 &  0   & -X_1 \\ 
\,X_3 & -X_2 &  X_1 &  0 
\end{pmatrix}\;. 
\end{gather*} 
\indent
Then, if $A={\rm Re}\,A+{\rm Im}\,A:S^3\to\Hq$ we have 
\begin{equation*} 
{\rm D}A=\dif^*({\rm Im}\,A)+\bigl(\dif({\rm Re}\,A)+(*\dif+2)({\rm Im}\,A)\bigr)\;. 
\end{equation*} 
Also, ${\rm D}^2=\D+2\,{\rm D}$ and the proof follows.  
\end{proof} 

\begin{rem} \label{rem:Beltrami} 
Proposition \ref{prop:solv} and its proof are similar to the discussion of the Beltrami fields  
equation on $\Re^3$ in \cite{KenPlu}\,. This discussion admits the following reformulation 
(which we do not imagine to be new):\\ 
\indent
Let $\D$ be the Hodge Laplacian on an oriented Riemannian 
three-manifold $N^3$\,. Then for a coclosed one-form $A$ the \emph{vector wave equation}  
$\D A=c^2A$, where $c$ is a constant, reads $$(*\dif+c)(*\dif-c)(A)=0\;.$$ 
Therefore \emph{we have a linear map from the space of coclosed solutions 
of the vector wave equation to the space of solutions of the} Beltrami fields equation  $$*\dif\!A=c\,A$$ 
\emph{given by} $$A\longmapsto (*\dif+c)(A)$$ 
\emph{which is easily seen to be surjective if} $c\neq0$\,.  
\end{rem} 

\indent
The following result is an immediate consequence of Proposition \ref{prop:solv}\,. 

\begin{cor} 
The set of global solutions to $\dif\!A+2*\!A=0$ on $S^3$ is equal to the space of left invariant 
one-forms. 
\end{cor} 

\begin{exm} 
If $A$ is a left-invariant one-form on $(S^3,h)$ considered with its canonical metric,  
then $g=\r^2h+\r^{-2}(\r\dif\!\r+A)^2$ is the Eguchi-Hanson metric II 
\cite{EgHa}\,. 
\end{exm}

\end{document}